\documentstyle[12pt,amscd]{amsart}
\newtheorem{Theorem}{Theorem}[section] \newtheorem{Lemma}[Theorem]{Lemma}
\newtheorem{Corollary}[Theorem]{Corollary}
\newtheorem{Proposition}[Theorem]{Proposition}

\def\To{\longrightarrow}
\def\reg{\operatorname{reg}}
\def\det{\operatorname{det}}

\def\pp{{\frak p}}
\def\M{{\cal M}}

\begin{document}
\title{On the asymptotic linearity \\ of Castelnuovo-Mumford regularity}
\author{Ng\^o Vi\^et Trung and Hsin-Ju Wang }
\address{Institute of Mathematics,  Box 631, B\`o H\^o, 10000 Hanoi, Vietnam}
\email{nvtrung@@thevinh.vnn.vn}
\address{Department of Mathematics,  National Chung Cheng University,
Chiayi 621, Taiwan}
\email{hjwang@@math.ccu.edu.tw}
\thanks{The first author is partially supported by the National Basic Research Program of Vietnam} 
\maketitle

\centerline{\small \it Dedicated to Wolmer Vasconcelos on the occasion of his sixtyfifth birthday} \vskip 0.7cm

\section*{Introduction} \smallskip

Let $A$ be a  commutative Noetherian ring with unity.
Let $R$ be a standard graded algebra over $A$, where ``standard"
means $R_0 = A$ and $R$ is generated by elements of $R_1$. Let
$R_+$ be the ideal generated by the elements of positive degree.
For any finitely generated graded $R$-module $M$ and $i \ge 0$ we
denote by
$H_{R_+}^i(M)$ the $i$th local cohomology module of $M$ with
respect to $R_+$. The {\it Castelnuovo-Mumford regularity} of $M$
is the invariant $$\reg(M) := \max\{a(H_{R_+}^i(M))+i|\ i \ge
0\},$$ where $a(H)$ denotes the maximal non-vanishing degree of a
graded $R$-module $H$ with the convention $a(H) = -\infty$ if $H =
0$. It is a natural extension of the usual definition of the 
Castelnuovo-Mumford regularity in the case $R$ is a graded algebra 
over a field.
 If $R$ is a polynomial ring over a field $k$, then
$M$ has a finite graded minimal free resolution: $0 \To F_r \To
\cdots \To F_1 \To F_0 \To M$ and $\reg(M) =
\max\{b_i(M)-i|\ i \ge 0\}$, where $b_i(M)$ denotes the maximal
degree of the generators of $F_i$. Hence one may view $\reg(M)$ as
a measure for the complexity of the structure of $M$ (see e.g.
\cite{BaM} and \cite{EG}). \par

There have been a surge of interest on the
behavior of the function $\reg(I^n)$ where $I$ is a homogeneous
ideal in a polynomial ring over a field. It was first discovered
by Bertram, Ein, and Lazarsfeld \cite{BEL} that if $I$ is the
defining ideal of a smooth complex variety,  $\reg(I^n)$ is
bounded by a linear function. Later, Swanson \cite{S} proved this
for all homogeneous ideals. Recently, Cutkosky et al \cite{CuHT}
and Kodiyalam \cite{K} found out that $\reg(I^n)$ is
asymptotically a linear function. See also \cite{Ch}, \cite{Cu},
\cite{CuEL}, \cite{GGP} for related results.
\par

The aim of this paper is to prove the above phenomenon for the
general case, namely, when $I$ is a graded ideal of  a standard
graded algebra $R$ over an arbitrary ring $A$. That is not a
simple task because the proof in the case $R$ is a polynomial ring
over a field is based strongly on the fact that finitely generated
modules over $R$ have finite minimal free resolutions,
whereas that fact does not hold in general. Moreover, as far as we
can see, the general case can not be reduced to this case. Our
approach  uses the characterization of the Castelnuovo-Mumford
regularity by means of filter-regular sequences combined with some
results originally due to \cite{CuHT} and \cite{K}. Even these
results need to be proved differently since the original proofs
are based on techniques which are not available in general such as
Gr\"obner basis and Nakayama's lemma. More general, we will
study the asymptotic regularity of the $I$-adic filtration of an arbitrary finitely generated graded $R$-module $M$. This generalization allows us to prove
that $\reg(\overline{I^n})$, where $\overline{I^n}$ denotes the
integral closure of $I^n$,  is asymptotically a linear function
with the {\it same slope} as $\reg(I^n)$ if $R$ is a domain over an affine ring. 
The latter fact was not known when $R$ is a
polynomial ring over a field. \par

We call a graded ideal $J \subseteq I$ an {\it $M$-reduction} of
$I$ if $I^{n+1}M = JI^nM$ for some $n \ge 0$. Define
$$\rho_M(I) := \min\{d(J)|\ \text{$J$ is an $M$-reduction of
$I$}\},$$ where for a graded module $N$, $d(N)$ denotes the minimum number $d$
such that there is a homogeneous generating set for $N$ with elements of degree
$\leq d$. Moreover, let
$\varepsilon(M)$ denote the smallest degree of the homogeneous elements of
$M$. Then our main result can be formulated as follows.\medskip

\noindent {\bf Theorem \ref{main}.}
{\it  Let $R$ be a standard graded ring over a  commutative Noetherian
ring with unity and $I$ a graded ideal of $R$. Let $M$ be a finitely
generated
graded $R$-module. Then there exists an integer $e \ge \varepsilon(M)$
such that for
all large $n$,} $$\reg(I^nM) = \rho_M(I)n + e.$$

\section{Filter-regular sequence and regularity} \smallskip

Let $A$ be an arbitrary commutative Noetherian ring with unity.
Let $R$ be a standard graded algebra over $A$ and $M$ a
finitely generated graded $R$-module.  \par

Let $z_1,\ldots,z_s$ be linear forms in $R$. We call $z_1,\ldots,z_s$
an {\it
$M$-filter-regular sequence} if $z_i \not\in P$ for any associated
prime $P \not\supseteq R_+$ of $(z_1,\ldots,z_{i-1})M$ for $i =
1,\ldots,s$. Filter-regular elements have their origin in the
theory of Buchsbaum rings \cite{SCT}.  \par

The Castelnuovo-Mumford regularity can be characterized by means of
filter-regular sequences. The following characterization was proved
implicitly in the proof of \cite[Theorem 2.4]{T2} for the case $M = R$.
But that proof also holds for any $M$.\par

\begin{Proposition} \label{regularity}
Let $z_1,\ldots,z_s$ be an $M$-filter-regular sequence of linear forms
which generate an $M$-reduction of $R_+$. Then
$$\reg(M) = \max\big\{a\big((z_1,\ldots,z_i)M:R_+/
(z_1,\ldots,z_i)M\big)|\ i = 0,\ldots,s\big\}.$$
\end{Proposition}

Note that a homogeneous ideal $Q \subseteq R_+$ is an $M$-reduction of
$R_+$ if and only if $(M/QM)_n = 0$  for all large $n$. \par

We shall see that any $M$-reduction of $R_+$ can be generated by
an $M$-filter-regular sequence in a flat extension of $A$.

\begin{Lemma} \label{flat}
Let $Q$ be an $M$-reduction generated by the linear forms
$x_1,\ldots,x_s$. For $i = 1,\ldots,s$ put $z_i = \sum_{j=1}^su_{ij}x_j$,
where $U = \{u_{ij}|\ i, j = 1,\ldots,s \}$ is a matrix of indeterminates.
Put
$$A' = A[U,\det(U)^{-1}],\ R' = R \otimes_AA',\ M' = M \otimes_AA'.$$
If we view $R'$ as a standard graded algebra over $A'$ and $M'$ as a
graded $R'$-module, then $z_1,\ldots,z_s$ is an $M'$-filter-regular
sequence.
\end{Lemma}

\begin{pf}
Since the elements $z_i$ are defined by independent sets of
indeterminates, it suffices to show that $z_1 \not\in P$ for any
associated prime $P \not\supseteq R_+'$ of $M'$. By the definition
of $R'$, such a prime $P$ must have the form  $\pp R'$ for some
associated prime $\pp \not\supseteq R_+$ of $M$. If $Q \subseteq
\pp$, then $(M/\pp M)_n = 0$ for all large $n$ since $M/\pp M$ is
a quotient module of $M/QM$. From this it follows that there is a
number $t$ such that $R_+^tM \subseteq \pp M$. Since
$\operatorname{ann}(M)\subseteq \pp$, this implies $R_+
\subseteq \pp$, a contradiction. So we
get $Q \not\subseteq \pp$. Since $Q = (x_1,\ldots,x_s)$, this
implies $z_1 = u_{11}x_1 + \cdots + u_{1s}x_s \not\in \pp R' = P$.
\end{pf}

Lemma \ref{flat} allows us to use Proposition \ref{regularity} for the
computation of $\reg(M)$. Indeed, since $A' = $ is a flat
extension of $A$, we have $H_{R_+'}^i(M')_n \cong
H_{R_+}^i(M)_n\otimes_AA'$ for all $n$ and $i \ge 0$,
whence
$$\reg(M) = \reg(M').$$

 From the characterization of the regularity by means of filter-regular sequences we can deduce the following relationship
between $d(M)$ and $\reg(M)$.
This was done in \cite[Proposition 4.1]{T1} for the case $R$ is a
polynomial ring over a local ring and $M$ is a graded ideal of
$R$. But the proof there also holds for the general case.

\begin{Proposition} \label{degree bound}
$d(M) \le \reg(M)$.
\end{Proposition}

\section{Bigraded module and regularity}

Let $S=A[X_1,\ldots,X_s,Y_1,\ldots,Y_v]$ be a polynomial ring over
a commutative Noetherian ring $A$ with unity. Given a sequence of
non-negative integers $d_1, \ldots,d_v$, we can view $S$ as a
bigraded ring with $\deg X_i = (1,0)$, $i = 1,\ldots,s$,  and
$\deg Y_j = (d_j,1)$, $j = 1,\ldots,v$. For convenience we will
assume that $d_v = \max\{d_i|\ i = 1,\ldots,v\}.$\par

Let $\M$ be a finitely generated bigraded module over
$A[X_1,\ldots,X_s,Y_1,\ldots,Y_v]$. For a fixed number $n$
put $$\M_n := \oplus_{a \ge 0}\M_{(a,n)}.$$
Then $\M_n$ is a finitely generated graded module over the
naturally graded polynomial ring $A[X_1,\ldots,X_s]$. We will
show that $\reg(\M_n)$ is asymptotically a linear function.
\par
If $s = 0$, then $\reg(\M_n)$ is the invariant
$$a(\M_n) = \max\{a|\ \M_{(a,n)} \neq 0\}.$$ This case is
settled by the following result which was already known in the
case $A$ is a field  \cite[ Theorem 3.3]{CuHT}. However, the
proof there used Gr\"obner basis (and integer programming) which
is not available in the general case.\par

\begin{Proposition} \label{s=0}
Let $\M$ be a finitely generated bigraded module over the
bigraded polynomial ring $A[Y_1,\ldots,Y_v]$. Then $a(\M_n)$ is
asymptotically a linear function with slope $\le d_v$.
\end{Proposition}

\begin{pf} The case $v = 0$ is trivial. Assume that $v \ge 1$.
Consider the exact sequence of bigraded $A[Y_1,\ldots,Y_v]$-modules:
$$0 \To [0_\M:Y_v]_{(a,n)} \To \M_{(a,n)} \overset{Y_v} \To
\M_{(a+d_v,n+1)} \To [\M/Y_v\M]_{(a+d_v,n+1)} \To 0.$$
Since $0_\M:Y_v$ and $\M/Y_v\M$ can be viewed as bigraded modules over
$A[Y_1,\ldots,Y_{v-1}]$,  using induction we may assume that
$a([0_\M:Y_v]_n)$ and $a([\M/Y_v\M]_n)$ are asymptotically
linear functions with slopes $\le d_v$.
As a consequence,
\begin{align*}
a([0_\M:Y_v]_n) + d_v & \ge a([0_\M:Y_v]_{n+1})\\
a(\M/Y_v\M)_n) + d_v & \ge  a([\M/Y_v\M]_{n+1})
\end{align*}
for all large $n$. Further, we are done if $a(\M_n)
= a([0_\M:Y_v]_n)$ for all large $n$.
\par
Since $a(\M_n) \ge a([0_\M:Y_v]_n)$ for all $n$, it remains to
consider the case that there exists an infinite sequence of
integers $m$ for which $a(\M_m) > a([0_\M:Y_v]_m)$. Putting this
condition into the above exact sequence we get $$a(\M_{m+1}) =
\max\{a(\M_m) + d_v, a([\M/Y_v\M]_{m+1})\}.$$ On the other
hand, we have $$a(\M_m)+d_v \ge a(\M/Y_v\M)_m + d_v \ge
a([\M/Y_v\M]_{m+1})$$ for $m$ large enough. So we can find infinitely many 
integers $m$
 such that $$a(\M_{m+1}) = a(\M_m) + d_v >
a([0_\M:Y_v]_m) + d_v \ge a([0_\M:Y_v]_{m+1}).$$ Fix such an integer $m$. 
Similarly, we
can prove that for all $n \ge m$, $a(\M_n) > a([0_\M:Y_v]_n)$
and therefore $a(\M_{n+1}) = a(\M_n) + d_v$. So $a(\M_n)$ is
asymptotically a linear function with slope $d_v$ in this case.
\end{pf}

Now we are going to consider the general case.

\begin{Theorem} \label{s>0}
Let $\M$ be a finitely generated bigraded module over the above bigraded
polynomial ring $A[X_1,\ldots,X_s,Y_1,\ldots,Y_v]$. Then $\reg(\M_n)$
is asymptotically a linear function with slope $\le d_v$.
\end{Theorem}

\begin{pf} By Proposition \ref{s=0} we may assume that $s \ge 1$.
For $i = 1,\ldots,s$ put $z_i = \sum_{j=1}^su_{ij}X_j$, where $U =
\{u_{ij}|\ i,j = 1,\ldots,s\}$ is a matrix  of indeterminates. Put
$A' = A[U,\det(U)^{-1}]$ and $\M' = \M \otimes_A A'.$ Then
$\M'$ is a finitely generated bigraded module over
$A'[X_1,\ldots,X_s,Y_1,\ldots,Y_v]$. Since $\M_n' = \M_n
\otimes_AA'$, we have $$\reg(\M_n) = \reg(\M_n')$$ for all $n
\ge 0$. By Lemma \ref{flat}, $z_1,\ldots,z_s$ is an
$\M_n'$-filter-regular sequence. Put $R = A'[X_1,\ldots,X_s]$.
Since $(z_1,\ldots,z_s) = (X_1,\ldots,X_s) = R_+$, applying
Proposition \ref{regularity} we get $$ \reg(\M_n') =
\max\big\{a\big((z_1,\ldots,z_i)\M_n':R_+/
(z_1,\ldots,z_i)\M_n'\big)|\ i = 0,\ldots,s\big\}. $$ It is easy
to check that $$ (z_1,\ldots,z_i)\M_n':R_+/(z_1,\ldots,z_i)\M_n'
= [(z_1,\ldots,z_i)\M':R_+/(z_1,\ldots,z_i)\M']_n. $$ Moreover,
we may consider $(z_1,\ldots,z_i)\M':R_+/(z_1,\ldots,z_i)\M'$ as
a graded module over the bigraded polynomial ring
$A'[Y_1,\ldots,Y_v]$. By Proposition \ref{s=0}, the invariant
$a\big([(z_1,\ldots,z_i)\M':R_+/(z_1,\ldots,z_i)\M']_n\big)$ is asymptotically a linear function with slope $\le d_v$ for $i =
0,\ldots,s$. Therefore, we can conclude that $\reg(\M_n)$ is
asymptotically a linear function with slope $\le d_v$.
\end{pf}

\section{Regularity of adic filtrations}

Let $A$ be a commutative Noetherian ring with unity.
Let $R$ be a standard graded algebra over $A$ and $I$ a graded
ideal of $R$.  We will apply the results
of the preceding sections to study the function $\reg(I^nM)$ for any
finitely generated graded $R$-module $M$.  \par

First we will establish a relationship between $d(I^nM)$ and the
invariants $\rho_M(I)$, $\varepsilon(M)$ introduced in the
introduction. This relationship was already known for the case $M = R$
and $R$ is a polynomial ring over a field \cite[Proposition
4]{K}. However the original proof does not hold in the general case
since it uses
Nakayama's lemma.

\begin{Lemma}\label{reduction degree}
 $d(I^nM)  \ge \rho_M(I)n+\varepsilon(M)$ for all $n \ge 0$.
\end{Lemma}

\begin{pf} We may write $I = J+K$, where $J$ and $K$ denote the ideals
generated by the homogeneous elements of $I$ of degree $< \rho_M(I)$
and $\geq  \rho_M(I)$, respectively. Then
$$I^nM=  JI^{n-1}M +K^nM.$$
Note that $K^nM$ is generated by homogeneous elements of degree
$\ge \rho_M(I)n+\varepsilon(M)$.
If $d(I^nM) < \rho_M(I)n+\varepsilon(M)$, then $I^nM=JI^{n-1}M$.
Hence $J$ is an $M$-reduction of $I$. Since by definition of $J$,
$d(J) < \rho_M(I)$,
this gives a contradiction to the definition of $\rho_M(I)$.
\end{pf}

Now we are going to prove the main result of this paper.

\begin{Theorem} \label{main}
Let $R$ be a standard graded ring over a  commutative Noetherian
ring with unity and $I$ a graded ideal of $R$. Let $M$ be a finitely
generated
graded $R$-module. Then there exists an integer $e \ge \varepsilon(M)$
such that for
all large $n$, $$\reg(I^nM) = \rho_M(I)n + e.$$
\end{Theorem}

\begin{pf}
Let $J$ be an $M$-reduction of $I$ with $d(J) = \rho_M(I)$.
Let $R[Jt] = \oplus_{n \ge 0}J^nt^n$ be the Rees algebra of $R$ with
respect to
$J$ and $\M = \oplus_{n \ge 0} I^nM$.
Since $I^{n+1}M = JI^nM$ for all large $n$, we may consider $\M$ as a
finitely generated graded module over $R[Jt]$. Assume that $R_1$ is
generated by the linear forms
$x_1,\ldots,x_s$ and $J$ is generated by the forms
$y_1,\ldots,y_v$ with $\deg y_j = d_j$, $j = 1,\ldots,v$.
Represent $R[Jt]$ as a quotient ring of the bigraded polynomial
ring $A[X_1,\ldots,X_s,Y_1,\ldots,Y_v]$ with $\deg X_i = (1,0)$,
$i = 1,\ldots,s$,  and $\deg Y_j = (d_j,1)$, $j = 1,\ldots,v$.
Then $\M$ is a finitely generated bigraded module over
$A[X,Y]$ with $\M_n = I^nM$.  By Theorem \ref{s>0},
$\reg(I^nM)$ is
asymptotically a linear function with slope $\le
\max\{d_1,\ldots,d_v\} = d(J)$. Let $dn+e$ be this linear
function. By Proposition \ref{degree bound} and Lemma
\ref{reduction degree} we have
$$\reg(I^nM) \ge d(I^nM) \ge \rho_M(I)n+\varepsilon(M)$$
for all $n$. Therefore, $d \ge \rho_M(I)$.
Since $d \le d(J) = \rho_M(I)$, we can conclude that $d =
\rho_M(I)$ and $e \ge \varepsilon(M)$. \end{pf}

Setting $M = R$ we immediately obtain the following consequence (see
\cite[Theorem 3.1]{CuHT} and \cite[Theorem 5]{K} for the case when $R$
is a polynomial ring over a field).

\begin{Corollary} Let $R$ be a standard graded ring over a
commutative Noetherian
ring with unity and $I$ a graded ideal of $R$.  Then there exists an
integer $e \ge 0$ such that for all large $n$,
$$\reg(I^n) = \rho_R(I)n + e.$$
\end{Corollary}

In particular, we may apply Theorem \ref{main} to study the
asymptotic regularity of the integral closures $\overline{I^n}$ of $I^n$ (see
\cite[Corollary 3.5]{CuHT} for the case $R$ is a polynomial ring
over a field where the slope could not be determined).

\begin{Corollary} \label{closure}
Let $R$ be a standard graded domain over over a
commutative Noetherian
ring with unity and $I$ a graded ideal of $R$ such that
 $\overline {I^{n+1}} = I \overline {I^n}$ 
for $n$ sufficiently large. Then $\reg(\overline{I^n})$ is asymptotically
a linear function with slope $\rho_R(I).$
\end{Corollary}

\begin{pf}
By assumption, there exists an integer $n_0$ such that
$\overline {I^{n+1}} = I \overline {I^n}$ for all $ n \ge n_0$. Put
$M = \overline{I^{n_0}}$.
Then $I^nM = \overline{I^{n+n_0}}$ for all $ n \ge 0$. As a consequence,
a graded ideal $J\subseteq I$ is an $M$-reduction if and only if
$\bar J = \bar I$ or, equivalently,  $J$ is an $R$-reduction of $I$. Therefore, $\rho_M(I) = \rho_R(I)$.
Hence the conclusion follows from Theorem \ref{main}.
\end{pf}

It is well-known that the assumption of Corollary \ref{closure} is satisfied for any ideal $I$ if $R$ is an affine domain.

\end{document}